\documentclass{mm}
\usepackage{amssymb, amsmath, amsfonts}
\usepackage[pdftex]{hyperref}
\usepackage[latin1]{inputenc}
\usepackage[fixlanguage]{babelbib}
\selectbiblanguage{spanish}
\usepackage{graphicx}
\usepackage{rotating}
\usepackage{url}
\usepackage{wrapfig}
\decimalpoint

\usepackage{lastpage}
\usepackage{enumerate}
\usepackage{paralist}

\begin{document}
\title{D\'e usted una buena pl\'atica de matem\'aticas}
\author[Daniel Pellicer]{Daniel Pellicer\\
\normalsize Centro de Ciencias Matemáticas\\
\normalsize UNAM - Morelia\\
\normalsize {\sf pellicer@matmor.unam.mx}
}

\begin{abstract}
\noindent Este texto busca brindar al lector un punto de partida para dar buenas pl\'aticas acerca de temas matem\'aticos.
\end{abstract}

\maketitle

\section{?`Est\'a dirigido a usted este texto?}

Este texto fue escrito pensando en ayudar a todo aquel que quiere mejorar su manera de dar pl\'aticas de matem\'aticas.

Le recomiendo no continuar leyendo si usted espera encontrar una receta fiel que garantice que su siguiente pl\'atica ser\'a un \'exito total. En cambio, si usted est\'a interesado en descubrir por usted mismo un camino que le permita mejorar sus exposiciones, el presente texto gustoso le ofrece un punto de partida. No espere encontrar una poci\'on m\'agica que solucione de golpe todas las deficiencias de sus charlas; en vez de eso esm\'erese en mejorarlas poco a poco hasta alcanzar un est\'andar con el que se sienta satisfecho. Este texto por s\'{\i} mismo no le har\'a dar mejores pl\'aticas; su esfuerzo y dedicaci\'on en mejorar s\'{\i} lo har\'a.



\section{?`Por qu\'e va usted a dar una pl\'atica?}

\begin{figure}[h]
\centering
\includegraphics[width=0.45\textwidth]{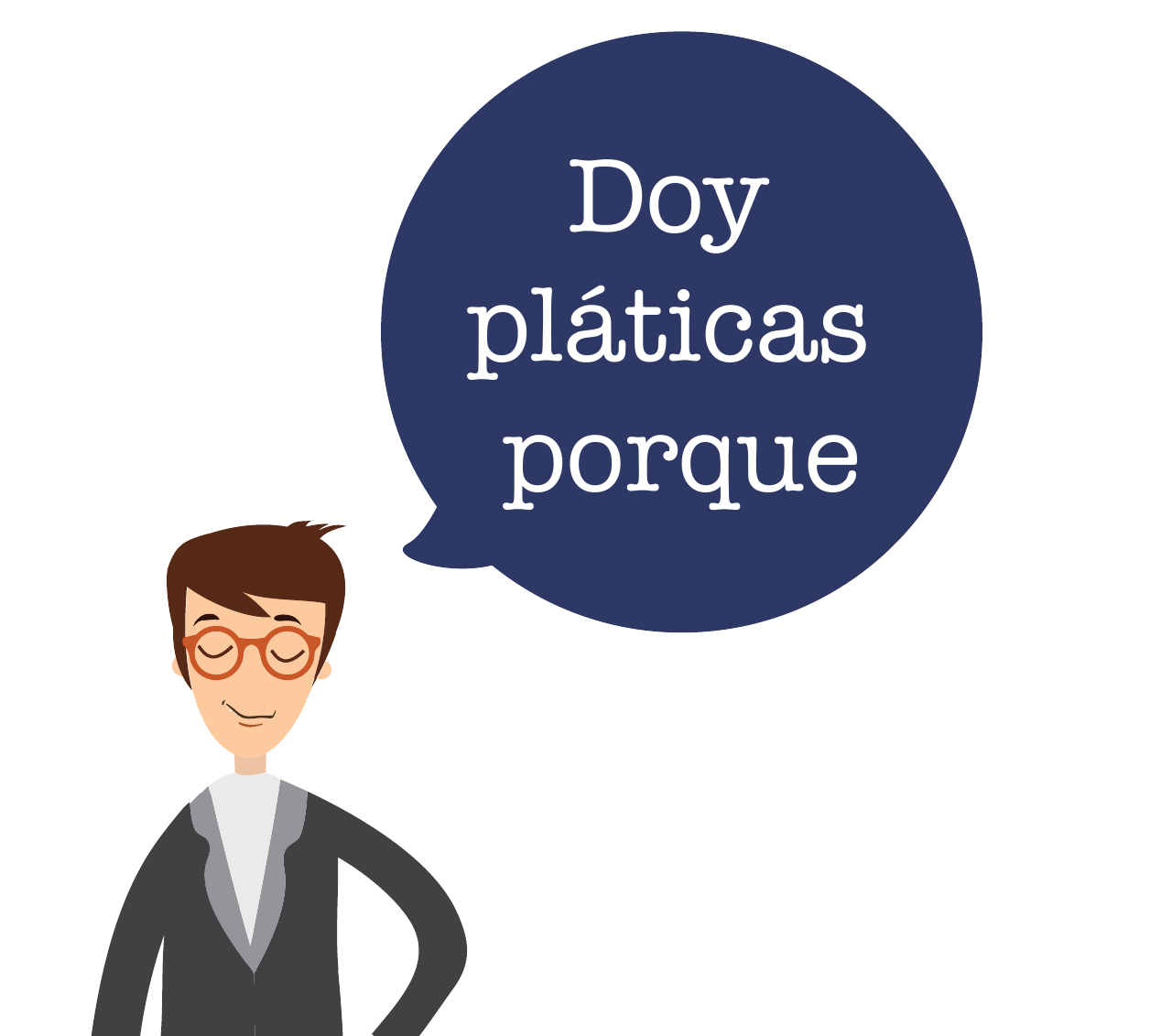}
 \end{figure}

Hay muchas razones por las que damos pl\'aticas de matem\'aticas, y el peso que le damos a estas razones var\'{\i}a de persona a persona. Le presento a continuaci\'on algunas de ellas.

{\em Doy pl\'aticas porque es un requisito.} Con frecuencia el programa acad\'emico en el que estamos inscritos o nuestro contrato la\-bo\-ral requiere que impartamos una o m\'as pl\'aticas.

{\em Doy pl\'aticas porque quiero transmitir mis conocimientos.} El gusto que tengo por el \'area de matem\'aticas que estudio me motiva a querer compartir mis conocimientos con mis compa\~neros estudiantes y profesores.

{\em Doy pl\'aticas para causar una buena impresi\'on.} Quiero mostrarle a expertos de mi \'area y a compa\~neros estudiantes que s\'e mucho del tema que estudio o que mis resultados tienen buena calidad.

{\em Doy pl\'aticas porque me siento comprometido a hacerlo.} Mi director de tesis (o alguna otra persona de renombre) me invit\'o a dar pl\'atica y no quiero rechazar su invitaci\'on.

{\em Doy pl\'aticas porque me gusta dar pl\'aticas.} Entre la atenci\'on que recibe el expositor, el estar frente a una audiencia y el platicar sobre mis temas favoritos de matem\'aticas, hay algo que me motiva a presentar ponencias.

Sin importar cu\'ales sean sus razones para dar pl\'aticas, debe tomarlas en cuenta al planear sus presentaciones. En secciones subsecuentes se detallan ejemplos frecuentes de maneras de sabotear las razones propias para dar pl\'aticas.

En este texto se da por hecho que los objetivos de su charla incluyen el lograr que los asistentes asimilen las ideas y resultados de los que va a hablar. Gran parte del contenido siguiente carece de sentido si esto no fuera relevante para usted cuando da pl\'aticas.

\section{?`Qui\'en va a escuchar la pl\'atica que usted dar\'a?}

Imag\'{\i}nese usted dando una pl\'atica a ni\~nos de primaria, y compare esa experiencia con una exposici\'on de su trabajo ante expertos en el \'area. Compare, por ejemplo, los temas de los que hablar\'{\i}a y el vocabulario que utilizar\'{\i}a ante tan distintas audiencias.

El ejercicio anterior le dejar\'a claro que distintos foros requieren distintas pl\'aticas, y por lo tanto, distintas preparaciones. Conf\'{\i}o en que usted se rehuse terminantemente a dar una misma pl\'atica en una escuela primaria y en un congreso de investigaci\'on.

Veamos las diferencias relevantes entre las dos audiencias descritas ante\-rior\-mente. En la primera los asistentes apenas tendr\'{\i}an noci\'on de los conceptos matem\'aticos m\'as b\'asicos, mientras que la segunda podr\'{\i}a contar con algunas de las personas con mayores conocimientos a nivel mundial en el tema a exponer. La primera no cuenta con conocimientos de \'areas afines que permitan relacionar el tema a abordar, mientras que la segunda ver\'a con agrado el uso de \'areas afines de matem\'aticas u otras ciencias para complementar la presentaci\'on. Los individuos de la primera audiencia no tienen madurez suficiente para comprender la pertinencia de la investigaci\'on realizada o en curso, mientras que la totalidad de la segunda est\'a directamente involucrada en el mundo acad\'emico.

\begin{wrapfigure}{r}{0.45\textwidth}
\centering
\includegraphics[width=0.45\textwidth]{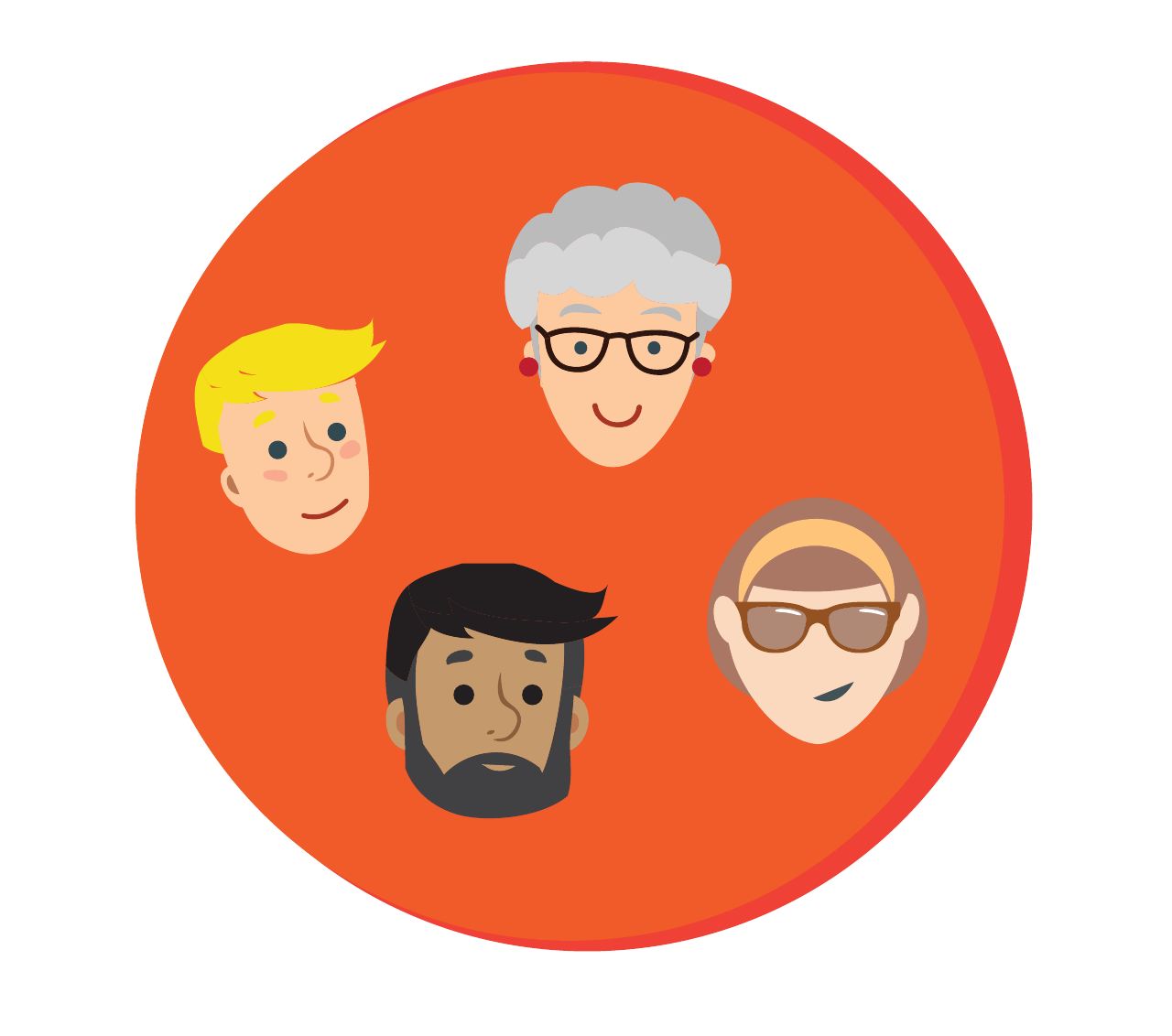}
 \end{wrapfigure}

Las dos audiencias reci\'en consideradas son solo casos extremos de un abanico de posibles p\'ublicos ante los cuales se presentan pl\'aticas de mate\-m\'a\-ticas semana a semana.
!`Pero no se requiere que sean tan radicalmente opuestas para que sea necesario preparar de manera distinta las pl\'aticas a presentar! Tenga en consideraci\'on que
\begin{itemize}
 \item la cantidad de conocimientos,
 \item la madurez,
 \item la rapidez para asimilar ideas nuevas y
 \item la comprensi\'on de la pertinencia de la investigaci\'on
\end{itemize}
son muy distintos entre un alumno a la mitad de la licenciatura, un alumno a finales de la maestr\'{\i}a y un investigador del \'area. Cada uno de ellos requiere una manera distinta de abordar el tema a exponer.


No siempre es posible dar gusto a todos los que escuchan nuestras pl\'aticas, en especial cuando sus niveles de conocimiento son muy variados. Es posible que o bien aquellos con m\'as conocimientos se aburran, o bien aquellos con menos conocimientos se pierdan. En ocasiones se puede librar esta dificultad haciendo una s\'{\i}ntesis de los conocimientos necesarios que cumpla las sigientes caracter\'{\i}sticas:
\begin{itemize}
 \item Debe ser suficientemente accesible para que quien no domina esos temas pueda desarrollar una intuici\'on de ellos. Ello requiere cierto ingenio. !`Recitar definiciones formales largas generalmente tiene el efecto de confundir en vez de crear intuici\'on!
 \item Debe incluir motivaci\'on y ejemplos interesantes, buscando que los expertos tambi\'en puedan tener inter\'es en seguir esta parte de
la charla.
 \item Las ideas deben seguir el flujo natural que requiere quien las escucha por primera vez. 
 \item No debe asumir que nuestra audiencia comprender\'a en un minuto algo que a nosotros nos tom\'o cinco minutos o m\'as en comprenderlo.
\end{itemize}
Indudablemente hacer lo anterior resta tiempo al tema del que en realidad queremos hablar. Hacerlo tampoco garantiza que quien no conoce los conceptos vaya a desarrollar suficiente intuici\'on para comprender lo que el expositor desea expresar en la pl\'atica, ni que quien ya conoce los conceptos no pierda inter\'es; sin embargo mejora la probabilidad de que m\'as asistentes comprendan y disfruten la pl\'atica.

En caso de que usted decida no esforzarse de m\'as para dejar contentos a todos los asistentes a su charla, y en vez de eso enfocarse solo a aquellos de un nivel acad\'emico determinado homog\'eneo, dicho nivel debe establecerlo de acuerdo al tipo de charla que le ha sido ofrecida. Es decir, si un experto en el \'area entra a una pl\'atica dirigida a estudiantes de licenciatura, usted debe asumir que el experto ser\'a paciente para permitir que los dem\'as asistentes conozcan los temas preliminares antes de llegar a la parte de inter\'es para el experto, aun si esto implica que solo los \'ultimos 5 minutos ser\'an relevantes para el experto. Por otro lado, en una pl\'atica de un seminario de investigaci\'on es de esperarse que usted pueda hablar de sus resultados m\'as avanzados. Un estudiante de licenciatura sin la preparaci\'on suficiente para seguir la pl\'atica deber\'a comprender que entr\'o al evento acad\'emico equivocado, o antes de estar debidamente preparado para ello.
Usted no debe sentirse mal de que su pl\'atica no haya cumplido las espectativas de una minor\'{\i}a de la audiencia, sobre todo cuando el evento acad\'emico no estaba particularmente orientado para esa minor\'{\i}a.

Como ver\'a el lector, preparar una pl\'atica de matem\'aticas para una audiencia homog\'enea requiere entender las caracter\'{\i}sticas espec\'{\i}ficas del sector acad\'emico al que ellos pertenecen. Por otro lado, una pl\'atica para una audiencia heterog\'enea requiere m\'as cuidado.


\section{?`Cu\'anto va a durar su pl\'atica?}

En la mayor\'{\i}a de los casos, cuando se nos invita a dar una pl\'atica se nos indica de cu\'anto tiempo dispondremos para ello, o al menos un tiempo estimado.

Una pl\'atica de 50 minutos permite plantear el contexto hist\'orico del problema y motivarlo adecuadamente. En este tipo de pl\'aticas es posible enunciar dos o m\'as resultados pertinentes al mismo conjunto de definiciones y resultados preliminares. Usualmente uno encuentra este tipo de pl\'aticas en seminarios, coloquios y conferencias plenarias de congresos.

En pl\'aticas de 20 minutos no es posible abarcar muchos resultados. En ellas se debe priorizar cubrir todo lo necesario para que se comprenda el resultado principal a exponer, por lo que debe contener un n\'umero relativamente peque\~no de definiciones y resultados preliminares. Los reportes de tesis en congresos y pl\'aticas de sesiones especiales suelen ser de esta duraci\'on aproximadamente.

\section{?`De qu\'e tema debe hablar?}

\begin{wrapfigure}{r}{0.45\textwidth}
\centering
\includegraphics[width=0.45\textwidth]{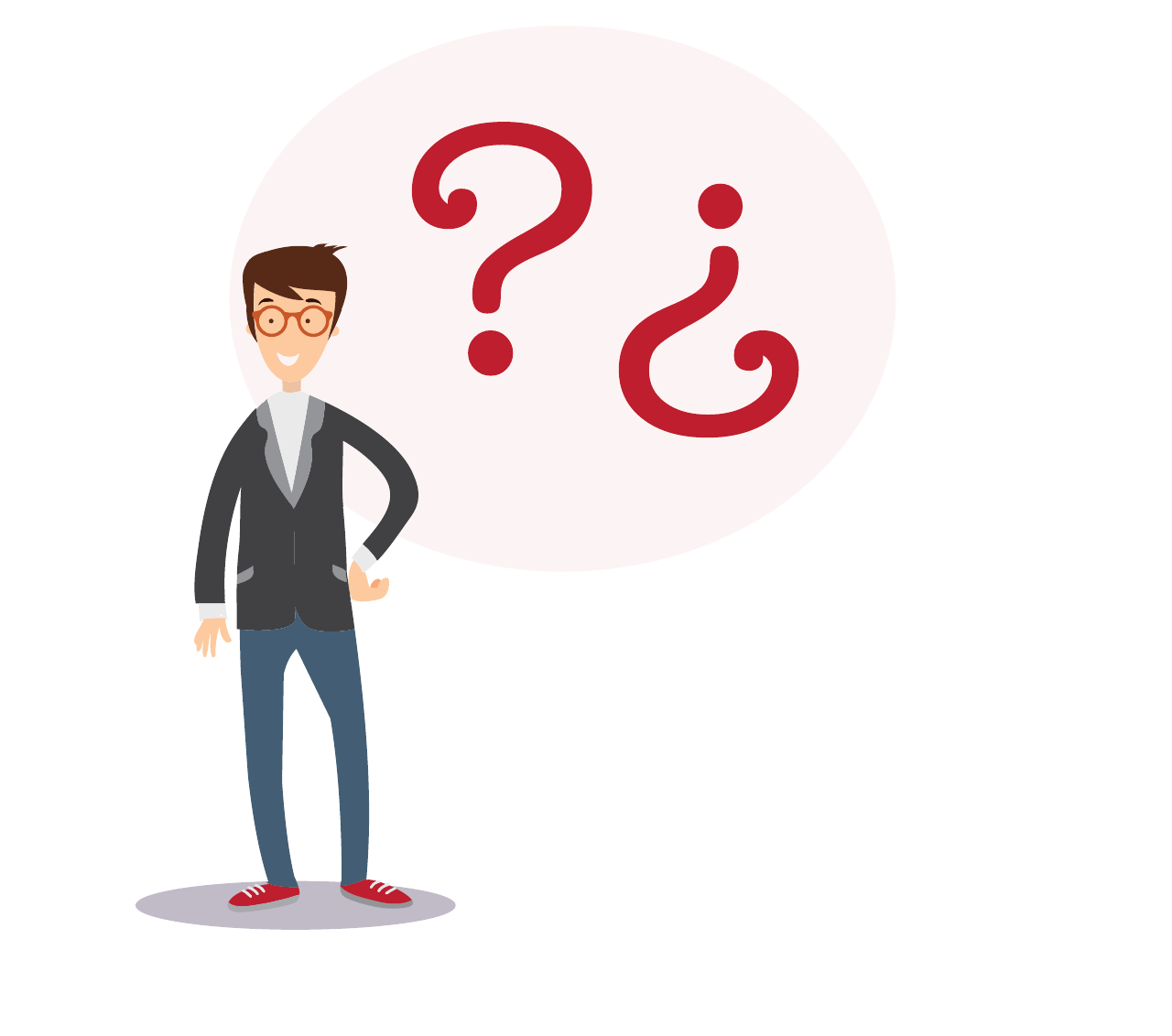}
 \end{wrapfigure}

Cuando se nos invita a dar pl\'aticas de matem\'aticas en congresos o seminarios, es conveniente conocer el tipo de audiencia esperada y la duraci\'on de la pl\'atica antes de decidir el tema del que hablaremos. Decidir primero el tema y despu\'es ajustarlo al tiempo disponible y al p\'ublico suele ocasionar que no se logren transmitir las ideas para las que la pl\'atica fue pensada.

Es natural tener la tentaci\'on de hablar acerca de nuestros logros m\'as recientes, o de nuestros resultados m\'as relevantes. Lo anterior no debe causar mayor dificultad en un evento dirigido a expertos en dichos temas. Sin embargo las dificultades se presentan de manera natural si el evento es dirigido a estudiantes o la audiencia incluye investigadores de otras \'areas de matem\'aticas.

Imagine primero que Andrew Wiles decidiera hablar de su demostraci\'on del \'ultimo teorema de Fermat. Este tema pudiera ser viable ante cualquier audiencia, aunque la exposici\'on debe variar significativamente de una au\-dien\-cia a otra. Si se presenta frente a un grupo de investigadores familiarizados con curvas el\'{\i}pticas (o con alg\'un otro elemento clave en su demostraci\'on), Andrew podr\'a mostrar detalles finos de la t\'ecnica que \'el utiliz\'o. Por otro lado, esa misma pl\'atica no ser\'{\i}a comprendida si es impartida como conferencia magistral de media hora ante estudiantes de licenciatura; ante esta \'ultima audiencia podr\'{\i}a enfatizar la relevancia hist\'orica del problema, plantearlo con toda precisi\'on e ilustrar ideas utilizadas en la demostraci\'on, sin formalizarlas ni profundizar en ellas. En realidad, ante estudiantes de 2o a\~no de licenciatura le tomar\'{\i}a mucho tiempo (tal vez el tiempo total destinado a la pl\'atica) definir curvas el\'{\i}pticas y crear intuici\'on suficiente de modo que los resultados que mencione tengan sentido para los asistentes. Ante una audiencia a\'un menos experimentada que la conformada por estudiantes de licenciatura, el principal objetivo de su pl\'atica deber\'{\i}a ser plantear el \'ultimo teorema de Fermat de modo que todo el p\'ublico entienda el enunciado, sin intentar esbozar una idea de demostraci\'on. Se puede tomar como objetivo dar a entender la importancia hist\'orica y la dificultad inherente a un problema con planteamiento tan sencillo.

El \'ultimo teorema de Fermat es un tema apto para audiencias variadas y distintas duraciones debido en parte a que es un tema famoso del que muchas personas, y en especial la gran mayor\'{\i}a de los matem\'aticos, han o\'{\i}do hablar con anterioridad. Sin embargo, la caracter\'{\i}stica principal que le permite ser presentado ante audiencias con pocos conocimientos matem\'aticos es que su planteamiento requiere muy pocos prerrequisitos. Para entender el enunciado es suficiente comprender conceptos que se adquieren en secundaria y preparatoria; aun frente a p\'ublico que no los maneje bien, pueden ser explicados en poco tiempo de manera intuitiva utilizando ejemplos.

Si, por otro lado, usted demostrara la hip\'otesis de Riemann y quisiera hablar de tan destacado logro, se topar\'a con el problema de que son pocas las audiencias preparadas para entender su enunciado. Estudiantes de primer a\~no de licenciatura apenas asimilan el concepto de funci\'on de manera formal, y tienen dificultades si el dominio y contradominio no est\'an dentro de los n\'umeros reales. A muchos de ellos les toma muchos meses comprender bien los n\'umeros complejos y manejar con soltura sus operaciones b\'asicas. Sin duda la comprensi\'on del significado de la funci\'on
\[\zeta(s)=\sum_{n=1}^{\infty}\frac{1}{n^s}\]
con dominio en los n\'umeros complejos ser\'a un reto de alto grado de dificultad para ellos; es probable que algunos pasen el tiempo restante de la charla tratando de entender qu\'e sentido tiene una suma de una infinidad de n\'umeros, o el significado de elevar un entero a una potencia imaginaria.  Al usted le tomar\'{\i}a mucho tiempo abordar convenientemente los ingredientes necesarios de modo que el enunciado sea comprendido. Probablemente la pl\'atica fuera m\'as exitosa si usted destinara su totalidad a la historia de la hip\'otesis de Riemann y a su relevancia en las matem\'aticas del \'ultimo siglo y medio, aun si su enunciado solo se esboza sin esperar que se comprenda a plenitud en qu\'e consiste.

En vista de lo anterior, es previsible que ante una audiencia conformada de personas sin estudios universitarios en matem\'aticas se logre una mejor comunicaci\'on entre el expositor y quienes lo escuchan si se elige como tema los n\'umeros complejos y no la hip\'otesis de Riemann en s\'{\i}. Tome usted en cuenta que en una pl\'atica acerca de los n\'umeros complejos es posible mencionar el contexto hist\'orico y la relevancia de la hip\'otesis de Riemann sin tener que enunciarla.

Es momento de regresar nuestra atenci\'on a las razones por las que usted dar\'a su siguiente pl\'atica. Si va a presentar una tesis ante un jurado para obtener un grado, su objetivo principal deber\'a ser el convencer a sus si\-no\-dales de sus conocimientos y destreza en el tema. Su presentaci\'on deber\'a tratar exactamente de el contenido de su tesis, aun si m\'as de la mitad de los asistentes no tienen estudios universitarios en matem\'aticas.

Por otro lado, si su objetivo principal es que los asistentes sigan las ideas que expone, inevitablemente deber\'a tomar en cuenta el tiempo disponible y la audiencia esperada antes de definir el tema; de otro modo corre el riesgo de ser usted mismo el primero en poner obst\'aculos al cumplimiento del objetivo planteado por usted para su pl\'atica.

\section{?`Cu\'ales preliminares deber\'a incluir en su pl\'atica?}

\begin{wrapfigure}{r}{0.45\textwidth}
\centering
\includegraphics[width=0.45\textwidth]{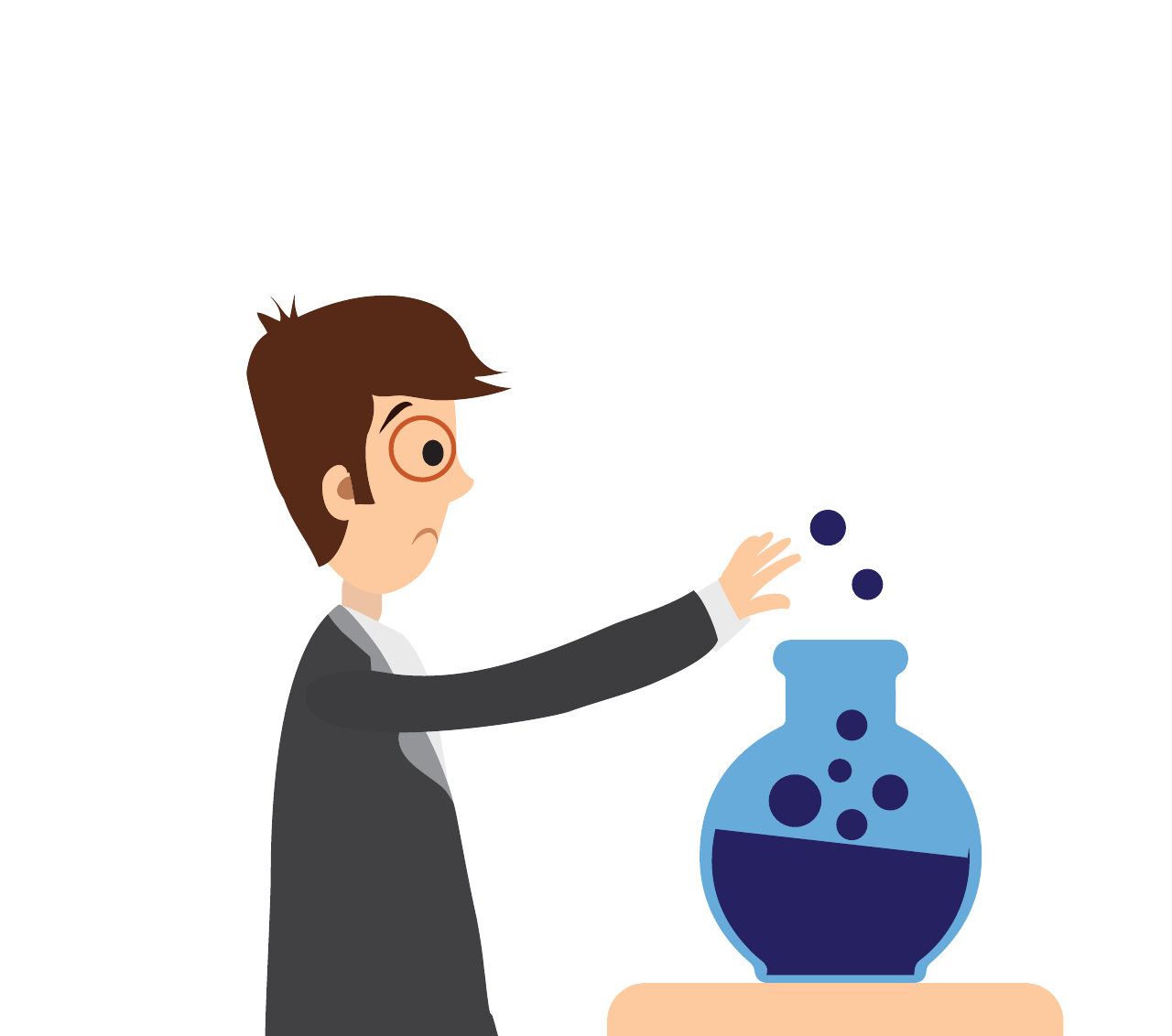}
 \end{wrapfigure}

Una vez establecido el tema de la pl\'atica debemos tener en mente el (los) resultado(s) o idea(s) a los que queremos llegar. El siguiente importante paso es decidir c\'omo vamos a llegar a ellos.

Si bien cuando escribimos un texto debemos buscar ser precisos y correctos, en un texto hablado no siempre es as\'{\i}. En un art\'{\i}culo de investigaci\'on los conceptos usados deben estar propiamente definidos, y debemos justificar nuestras afirmaciones con argumentos, demostraciones o referencias. En algunas pl\'aticas se espera un car\'acter similar al de un art\'{\i}culo de investigaci\'on, mientras que en otras se entiende que el tiempo y el nivel de conocimientos de la audiencia hacen que la rigurosidad sea inadecuada.

Imagine usted que lo invitan a hablar acerca de su tesis ante un grupo de expertos en el tema. En esa circunstancia es innecesario que defina algunos conceptos que se incluir\'{\i}an en art\'{\i}culos de investigaci\'on, dado que todos los asistentes a la pl\'atica los dominan. El contenido de dicha pl\'atica debe enfocarse en lo que es novedoso para esa audiencia, que puede incluir definiciones recientes, los resultados de su investigaci\'on, y detalles de las pruebas de dichos resultados.

Por otro lado, si va a hablar acerca de su tesis en un seminario de estudiantes, al que asisten estudiantes de \'areas distintas a la de usted, deber\'a incluir en la pl\'atica definiciones y/o motivaci\'on para los conceptos que no se estudian en el tronco com\'un de las licenciaturas y maestr\'{\i}as en matem\'aticas. Asumir que los asistentes saben el significado de homolog\'{\i}a, forcing, matroide o superficie de Riemann (por poner unos pocos ejemplos) y no definir o motivar estos conceptos adecuadamente suele ocasionar que los asistentes no puedan seguir la pl\'atica a partir de que estos conceptos se convierten en parte importante de esta. Lo anterior fomenta que en el tiempo restante cambien su atenci\'on a sus propios problemas matem\'aticos, a lo que har\'an despu\'es de la pl\'atica, o a sus dispositivos m\'oviles. Dado que tal pl\'atica comenzar\'{\i}a en un nivel m\'as bajo que la pl\'atica que se impartir\'{\i}a ante expertos, es de esperarse que no le alcance el mismo tiempo para cubrir el mismo material. Una alternativa es que priorice la comprensi\'on de los enunciados de sus resultados, dejando de lado las demostraciones de los mismos. Al hacerlo sacrificar\'a el intento de convencer a los asistentes de que su resultado es cierto, en favor de que asimilen lo que plantea su resultado, junto con la relevancia de este.

Una vez que usted ha decidido qu\'e definiciones incluir\'a en la pl\'atica, es importante distinguir entre aquellas que basta con enunciar y aquellas en las que debe abundar. Por ejemplo, si la audiencia consiste de estudiantes de segundo semestre de licenciatura, es posible definir de manera clara el complemento ortogonal de un subespacio vectorial de $\mathbb{R}^3$ usando para ello unos cuantos segundos. Por otro lado, si se va a hablar de la forma can\'onica de Jordan de una matriz cuadrada, ello requerir\'a, adem\'as de una definici\'on, suficiente motivaci\'on y ejemplos de modo que los resultados que mencionan dicho concepto tengan sentido para los asistentes.

En ocasiones un resultado que queremos presentar requiere de muchas definiciones y resultados previos, pero se simplifica considerablemente si lo explicamos con un caso particular en el que nuestro resultado no se trivialice. Dedicando la pl\'atica al caso particular damos m\'as oportunidad a quienes nos escuchan de que sigan las ideas claramente, y no excluimos la posibilidad de que al final se mencione el resultado general y se expliquen las diferencias. En cambio, intentar el caso general desde un principio puede ser demasiado para una audiencia poco conocedora del tema, y puede propiciar que se deje de seguir la pl\'atica antes de que se est\'e cerca de mencionar el resultado principal. Siempre contemple lo anterior en caso en que se le solicite hablar acerca de un tema espec\'{\i}fico ante una audiencia poco preparada para entenderlo.

Las cosas que son importantes o interesantes para usted no necesariamente lo son para la audiencia. Dar un amplio rodeo incluyendo temas no relacionados con la charla para justificar que cierto resultado es muy importante o muy interesante puede tener como consecuencia que los asistentes pierdan inter\'es. En particular, si da usted una pl\'atica de cierta \'area ante una audiencia con formaci\'on predominante en un \'area distinta, ponga entusiasmo en abundar en los temas de inter\'es com\'un (sin omitir las definiciones y motivaciones necesarias). Suele dar peores resultados pensar en `educar' a esa audiencia para que vean la belleza de lo que es atractivo para usted, si para ello debe tomar mucho tiempo, incluir muchas definiciones o hablar de varios detalles t\'ecnicos, distrayendo a la audiencia del objetivo real de la charla.

Tome por ejemplo un matem\'atico que recientemente obtuvo resultados en ecuaciones diferenciales con muchas aplicaciones en f\'{\i}sica de fluidos, y que gracias a ello es invitado a hablar en un congreso de dicha \'area de la f\'{\i}sica. Naturalmente en esas circunstancias se espera que la pl\'atica gire alrededor de la relevancia de los resultados, su contexto hist\'orico y sus potenciales aplicaciones, y a ello deber\'a dedic\'arsele el cuerpo de la exposici\'on. Es de esperar que la audiencia no vea bien si usted dedica la mayor parte de la pl\'atica a la demostraci\'on de los resultados, o en fundamentos de las ecuaciones dife\-ren\-cia\-les que no sean indispensables para la comprensi\'on de los enunciados demostrados. El mismo fen\'omeno, si bien en menor grado, ocurre cuando a un matem\'atico especializado en cierta \'area se le invita a dar una pl\'atica en un evento especializado en otra \'area.

Por \'ultimo, tome usted en cuenta que las personas a quienes les habla quieren escuchar acerca de las matem\'aticas que usted desarrolla, y no de las cuentas que fueron necesarias para llegar a sus resultados. Mostrar largas listas de igualdades que se sigan una de la otra suele provocar desinter\'es en entenderles a fondo, y en el contexto de una pl\'atica con frecuencia se pueden sustituir por explicaciones intuitivas de aquellas razones que est\'an detr\'as de la veracidad del resultado.

\section{Acerca del t\'{\i}tulo y resumen de su pl\'atica}

La elecci\'on del t\'{\i}tulo de nuestra ponencia no es un detalle menor, dado que puede influir en el mucho o poco inter\'es que la gente tenga en asistir.

Idealmente el t\'{\i}tulo de nuestra charla deber\'{\i}a cumplir simult\'aneamente las siguientes caracter\'{\i}sticas.
\begin{itemize}
 \item Ser breve.
 \item Dar una primera idea de lo que se debe esperar de la pl\'atica (por ejemplo, especificando el tema de la pl\'atica o sugiriendo el tono que va a tener).
 \item Invitar a todo aquel que pudiera interesarle la pl\'atica.
\end{itemize}
No siempre es sencillo cubrir todos estos aspectos, y debe usted tener en cuenta cu\'ales priorizar.

Hay varias razones por las cuales gente que usted quiere que asista a su pl\'atica puede decidir no ir. En congresos grandes en los que hay sesiones simult\'aneas los asistentes deben decidir a lo m\'as una de las varias pl\'aticas que se ofrecen en un mismo momento. En eventos m\'as peque\~nos puede ser que algunos asistentes acuerden darse tiempo para trabajar en proyectos pendientes, y para ello decidan ausentarse de unas pocas pl\'aticas.
A aquellas personas que tengan razones para ausentarse de algunas pl\'aticas y quisieran decidir con precauci\'on a cu\'ales asistir, usted les hace sencillo el trabajo si en su t\'{\i}tulo no pueden ver de qu\'e trata su pl\'atica. Esto puede deberse a que en el t\'{\i}tulo incluya muchos t\'erminos t\'ecnicos, o porque asuma de los asistentes conocimientos previos de muchos conceptos. Si no le es posible expresar todo lo que quer\'{\i}a en pocas palabras, !`deje esa parte para el resumen!

Como ejemplo, imagine que usted dar\'a una pl\'atica introductoria de teor\'{\i}a de m\'odulos, donde su principal objetivo es que otras personas, en particular estudiantes j\'ovenes, conozcan su \'area de estudio. Para ello piensa mostrar a detalle una amplia gama de ejemplos. Un t\'{\i}tulo como `M\'odulos simples, semisimples, inyectivos y proyectivos' env\'{\i}a el mensaje que la pl\'atica no ser\'a apta para quienes no conozcan ya esos conceptos (al menos algunos de ellos). En lugar de eso puede llamarla `Introducci\'on a los m\'odulos', `La diversidad de los m\'odulos', o de alguna otra manera (tan formal o informal como desee) que sugiera que todo estudiante de licenciatura es bienvenido. Por otro lado, si su charla es acerca de su tesis y quiere mostrar sus resultados, pero debido al tiempo disponible debe asumir que la gente est\'a familiarizada con m\'odulos, debe adoptar un t\'{\i}tulo como el primero, y no como los segundos.

Es deseable que el resumen cumpla las siguientes caracter\'{\i}sticas:
\begin{itemize}
 \item Describir el contenido de la pl\'atica. Si en el momento en que debe enviar el resumen aun no conoce con precisi\'on el contenido de su pl\'atica, puede enviar un resumen vago que incluya frases como `...se mencionar\'an resultados acer\-ca de...' sin especificar de qu\'e resultados se trata.
 \item No incluir cosas que no se van a abordar.
 \item Seguir un estilo similar al de la pl\'atica. Si en la pl\'atica va a avanzar poco a poco en los conceptos, en el resumen se deben incluir frases como `se abordar\'an las definiciones y propiedades b\'asicas de ...', y se pueden incluir algunas pocas nociones que ayuden a entender el tipo de matem\'aticas que se desarrollar\'an. Por otro lado si se asumen prerrequisitos, estos se deben sugerir en el resumen, al utilizar dichos t\'erminos sin definici\'on previa.
 \item Buena ortograf\'{\i}a y redacci\'on. Res\'umenes que no cumplen con este requisito suelen causar una mala imagen del autor. Esto puede predisponer a potenciales asistentes para decidir no acudir a la charla, o no darle la seriedad debida.
\end{itemize}



Un resumen no es sustituto de tiempo que se pueda ahorrar en la pl\'atica. No debe asumir que a la hora de su charla los asistentes recuerden lo que dice su resumen. Si las definiciones precisas ocupan mucho espacio e insiste en incluirlas en el resumen, busque descripciones intuitivas o vagas que den suficiente idea, y reserve la definici\'on formal para su charla.

No hay formatos est\'andares de resumen que sirvan para todas las charlas. Cada exposior debe elaborar el suyo propio antes de su pl\'atica. Puede guiarse en res\'umenes de eventos anteriores de car\'acter similar para ver una muestra de estilos y longitudes.

\section{Acerca de la cantidad de material en su pl\'atica}

Imag\'{\i}nese usted que lee una novela de 100 p\'aginas en la que aparecen 200 personajes con nombre y apellido. Cada uno de estos personajes se presenta con pocas situaciones a su alrededor, pues despu\'es de todo la novela tiene solo 100 p\'aginas. Sin duda usted tendr\'a problemas en memorizar qui\'en es qui\'en. Usted podr\'a leer la novela de corrido esperando poder detectar del contexto a qu\'e persona corresponde cada nombre, en ocasiones teniendo que superar la confusi\'on originada por rebautizar a alg\'un personaje err\'oneamente. Otros dos resultados probables en esta situaci\'on son: que usted pierda inter\'es en la novela y deje de leerla, o bien que con frecuencia regrese a p\'aginas anteriores para recordar qui\'en era la persona a quien pertenec\'{\i}a tal nombre, y en qu\'e contexto aparec\'{\i}a. Aun en este \'ultimo caso, es muy probable que usted se quede con la impresi\'on de que la trama del libro pudo haberse expuesto con menos nombres de personajes (tal vez incluso con menos personajes).

El mismo fen\'omeno ocurre en pl\'aticas, sin importar si son de 20, 30 \'o 50 minutos. En ellas los conceptos (actores de la trama) en ocasiones son introducidos por definiciones formales, en ocasiones por definiciones intuitivas, y en ocasiones el expositor espera que el p\'ublico los conozca de pl\'aticas o cursos anteriores. Piense que quiere evitar la analog\'{\i}a con la novela, y por ello deber\'a buscar que o bien en su pl\'atica haya pocas definiciones, o bien que cada definici\'on est\'e debidamente ejemplificada y motivada para causar la menor confusi\'on posible en quienes escuchan.

No hay receta universal para solucionar este problema cuando se presenta. Ejemplificar y motivar definiciones toma tiempo, y no debe exagerarse en ello si se quiere dar el debido tiempo al contenido principal de la pl\'atica.
Por otro lado, decidirse por incluir pocas definiciones puede ocasionar que no se entienda cuando menciona conceptos que no defini\'o. Por supuesto tampoco debe suponer que los asistentes est\'en familiarizados con conceptos de los que hasta ese momento poco o nada hayan escuchado.
Usted puede considerar las siguientes alternativas.

Suponga que algunas hip\'otesis de alg\'un teorema, si bien necesarias, no entran dentro del marco general del tema que se est\'a exponiendo, y son originadas por c\'alculos o por resultados previos de sabor distinto al de la pl\'atica. El expositor debe procurar invertir una porci\'on lo m\'as peque\~na posible del tiempo disponible para hacer menci\'on de dichas hip\'otesis. Una manera de hacerlo es, !`no enunciarlas! y sustituirlas por frases como `algunas propiedades t\'ecnicas'. Si alguno de los asistentes est\'a particularmente interesado en ellas podr\'a preguntar al respecto al final de la charla, momento en el cual es posible extenderse en estas hip\'otesis sin interrumpir el flujo natural de las ideas planeadas.

Otra manera com\'un de evitar dar exceso de definiciones a memorizar a quien le escucha consiste en evitar algunos nombres que de inicio no dan informaci\'on a quien no los conoce (piense en la primera vez que escuch\'o `isomorfo', `inyectivo', `ortogonal'). En vez de ellos se puede usar calificativos informales que permitan al lector entender la idea que se quiere expresar. Por ejemplo, suponga que quiere dar una pl\'atica para estudiantes de nuevo ingreso a la licenciatura, y una de las hip\'otesis que requiere para su resultado principal es que cierto par de vectores sean linealmente independientes, pues de lo contrario aparecer\'a una divisi\'on entre cero. Puede ejemplificar lo que ocurrir\'{\i}a si los vectores fueran linealmente dependientes (sin usar esos t\'erminos) y concluir que al conjunto de vectores deseados se les llamar\'a `buenos', o alternativamente se dir\'a que son `malos' los linealmente dependientes. Busque dejar claro que tales t\'erminos se utilizar\'an \'unicamente durante la charla, y quienes quieran profundizar en el tema deber\'an adoptar los t\'erminos convencionales.

En textos escritos es com\'un agrupar las definiciones al inicio y despu\'es usar los conceptos seg\'un se van necesitando. En un texto hablado no hay oportunidad de `regresar unas p\'aginas' para leer la definici\'on de la que se nos est\'a hablando y que ya hemos olvidado; es m\'as conveniente dar y motivar cada definici\'on cuando se le va a usar.

De igual manera que se le invit\'o a reducir el n\'umero de definiciones a un nivel que el p\'ublico pueda asimilar f\'acilmente, se hace una invitaci\'on an\'aloga para el n\'umero de resultados. Algunos resultados preliminares no aportan ideas relacionadas al tema de la pl\'atica, y por lo mismo podr\'{\i}an ser omitidos, en todo caso haciendo menci\'on de que su existencia permite enunciar los resultados centrales. Frases como `resultados t\'ecnicos permiten concluir que...' son aceptables en una pl\'atica, y pueden ser aclarados posteriormente ante aquellos que muestren inter\'es en ello.

Por \'ultimo, se invita a usted a no caer en la tentaci\'on de hablar m\'as r\'apido buscando de esa manera cubrir m\'as material. Aumentar la velocidad del discurso tiende a generar consecuencias negativas ocasionadas por otorgar un menor tiempo a los asistentes para comprender las ideas que se les exponen. Piense que al incluir mayor material del debido corre el riesgo de que se comprendan menos ideas que si se hubiera expuesto una menor cantidad de material. Si es as\'{\i}, ?`tiene caso incluir ese material adicional?

\section{Preparaci\'on de su pl\'atica:}

Cuando usted ya ha decidido el tema de la pl\'atica y el car\'acter que se le va a dar, as\'{\i} como las definiciones y resultados que quiere abordar, puede proceder a elaborar la secuencia de ideas a expresar en la pl\'atica. Durante este proceso debe cuidar que toda menci\'on de un objeto matem\'atico hecha antes de que dicho objeto sea definido propiamente, debe tener un car\'acter informativo de hacia d\'onde se dirige la charla, y abordar el objeto \'unicamente de manera intuitiva. De igual manera, si piensa enunciar cierto resultado antes de haber descrito adecuadamente las componentes del resultado, h\'agale saber a quienes lo escuchan que esas componentes ser\'an explicadas m\'as adelante. Cuide de no dar el mensaje de que no pusieron atenci\'on en definiciones cuando (hipot\'eticamente) las mencion\'o.

Le sugiero hacer el ejercicio de ponerse en el papel del p\'ublico. Revise el orden de las ideas que piensa expresar, intentando prever el grado en el que los asistentes pueden comprender cada idea dados sus conocimientos y lo que usted habr\'a dicho con anterioridad en la pl\'atica. Es com\'un darse cuenta de que ser\'{\i}a bueno cambiar el orden de las definiciones, o incluir motivaciones.

A veces nos sentimos tentados a incluir en la pl\'atica adelantos de c\'omo va a ser usado lo que acabamos de explicar. En ocasiones eso justifica adecuadamente una definici\'on rebuscada o un resultado t\'ecnico. Sin embargo, un adelanto as\'{\i} no contribuye a la pl\'atica si parte de lo dicho no es comprendido por su audiencia, ya sea porque va a ser definido m\'as adelante o porque no se ha introducido el tema debidamente. En esos casos es mejor omitir el adelanto.

Parte importante de la preparaci\'on de una pl\'atica es hacer un ensayo general, midiendo el tiempo que tarda en concluir su exposici\'on. Esto puede ser ante su tutor, en un seminario, o present\'andose la charla a usted mismo.
Es importante que usted le de importancia al tiempo que le fue asignado, y est\'e preparado a terminar sin hablar de todo lo que ten\'{\i}a pensado porque se agot\'o el tiempo, o a agregar detalles o ejemplos si el tiempo lo permite.
No olvide usted tomar en cuenta tiempo para preguntas, ya sea que sean formuladas en el transcurso la pl\'atica, o al final de ella.

Es de mal gusto terminar una pl\'atica con demasiada anticipaci\'on, pues deja la impresi\'on de que el tiempo disponible fue mal empleado.

Terminar una pl\'atica despu\'es del tiempo previsto es causa de problemas m\'as concretos. En ocasiones algunos asistentes tienen compromisos inmediatamente despu\'es de la hora dedicada a un seminario; extenderse en este tipo de pl\'aticas inevitablemente provocar\'a que dichas personas tengan que abandonar el sal\'on, causando distracci\'on en la audiencia y el expositor. Es importante respetar el tiempo de los asistentes tambi\'en en eventos que tienen sesiones simult\'aneas; si alguien plane\'o escuchar una pl\'atica posterior a la de usted que se desarrollar\'a en otra sesi\'on, a esta persona no le ser\'a posible asistir a ambas pl\'aticas en caso de que usted tarde en terminar su charla. Finalmente, es aconsejable no crear dificultades a los organizadores, por ejemplo en caso de que el recinto en el que se lleva a cabo el evento acad\'emico deba ser abandonado puntualmente a la hora programada de t\'ermino de la \'ultima pl\'atica (sea la de usted, o la de alguien que hablar\'a despu\'es de usted).



En caso de que termine el ensayo de su pl\'atica varios minutos antes de cumplirse el tiempo disponible, le sugiero elaborar una lista de ejemplos, resultados o precisiones que sean candidatos a ser agregados; estos deber\'an ser priorizados de acuerdo tanto a la relevancia con el tema general de la charla, como a lo \'util que ser\'a para facilitar la comprensi\'on del material, e incluir \'unicamente aquellos que est\'en m\'as altos en la lista de prioridades, de modo que no se exceda el tiempo de la pl\'atica. Evite la tentaci\'on de incluir una idea que aleje a quienes le escuchan del tema alrededor del que va a hablar, sobre todo si esa idea requiere agregar definiciones y resultados preliminares innecesarios.

El caso contrario es m\'as frecuente. Ocurre seguido que al practicar nuestra pl\'atica nos encontramos con que rebasamos el tiempo ofrecido. En esas ocasiones es necesario eliminar partes de la pl\'atica, y este proceso requiere cierto cuidado. Es conveniente identificar primero si hemos incluido definiciones o resultados superfluos que distraigan la atenci\'on hacia temas distintos al tema principal a tratar; estos deber\'an ser los primeros candidatos a ser eliminados del contenido. Posteriormente se pueden buscar dentro del contenido de la pl\'atica definiciones y resultados que sean demasiado t\'ecnicos para que la explicaci\'on que hacemos de ellos aporte sustantivamente a la presentaci\'on; dicha explicaci\'on puede ser sustituida por una breve menci\'on, aclarando que se omiten los detalles debido a su car\'acter t\'ecnico. Si lo anterior a\'un no es suficiente y la pl\'atica todav\'{\i}a resulta demasiado larga, se puede reducir el n\'umero de definiciones formales, intercambiando algunas de ellas por ideas intuitivas; hay que elegir particularmente aquellos conceptos cuya definici\'on formal requiera mucha explicaci\'on que no est\'a directamente ligada a la idea general a tratar, y su esencia pueda ser transmitida por una idea intuitiva que requiere poco tiempo para describirla. En caso de que despu\'es de seguir las sugerencias anteriores la exposici\'on a\'un tome dema\-sia\-do tiempo deber\'a usted considerar seriamente en reducir el contenido de la pl\'atica, eliminando los resultados o ideas que aporten menos al tema a tratar. Alternativamente, en ocasiones se puede exponer el mismo contenido, pero sin la generalidad planeada originalmente; se pueden seguir las definiciones y resultados para un caso particular, y solo al final de la pl\'atica mencionar direcciones en que se generaliza. Esto resulta \'util cuando las definiciones son complicadas y requiren muchos ejemplos y motivaciones, pero existen casos particulares sencillos o bien conocidos en los que se pueden ilustrar los resultados que se quieren comunicar.

Aun cuando usted haya medido su pl\'atica y la duraci\'on sea \'optima deber\'a tener en cuenta que algunos imprevistos pueden forzarle a terminar su pl\'atica antes de cubrir todo lo previsto (muchas preguntas durante la exposici\'on, retraso al iniciar, etc). Busque que, de encontrarse en la situaci\'on de que le quedan 5 minutos y todav\'{\i}a falta mucho por exponer, sea sencillo decidir qu\'e partes omitir de modo que no deje fuera las ideas m\'as importantes.

\begin{figure}[h]
\centering
\includegraphics[width=8cm]{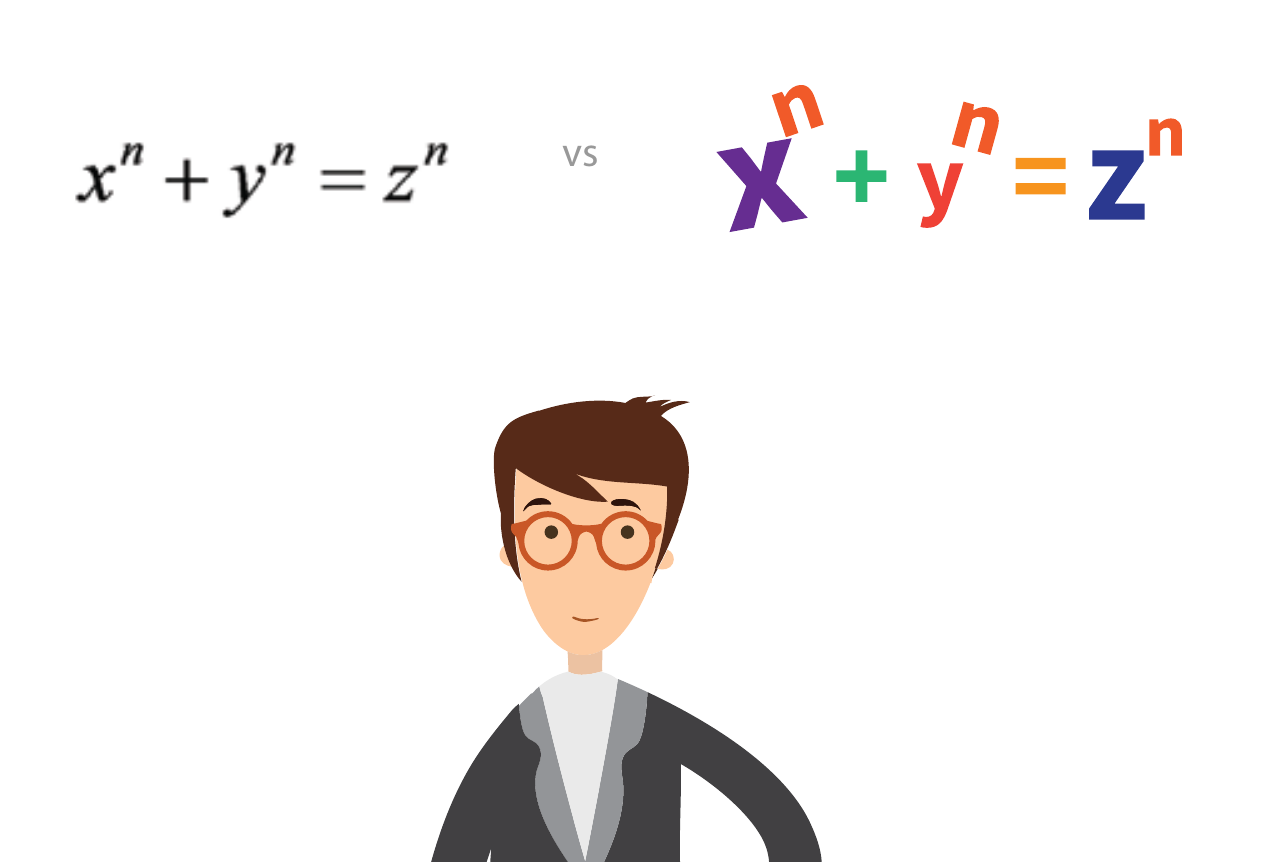}
 \end{figure}

\section{Elecci\'on de palabras que usar\'a en su pl\'atica}

Cada uno de nosotros tiene su propio car\'acter y sus propios gustos. Estos son un factor determinante en el estilo que elegimos para presentar nuestros resultados, siempre que dicha elecci\'on no est\'e supeditada a requerimientos del foro en el que hablaremos.

Si usted tiene facilidad para bromear entre sus conocidos, puede intentar incluir chistes y comentarios graciosos en sus pl\'aticas; esto puede tener como consecuencia que tanto usted como su audiencia se relajen del rigor del tema expuesto y la presentaci\'on resulte m\'as placentera.

En cambio, si bromear no es lo suyo no tiene por qu\'e forzar chistes en sus presentaciones, y en cambio puede dar una presentaci\'on tan sobria y formal como lo desee. Hacer algo a lo que no est\'a habituado y que no le entusiasma particularmente puede cargarle con una preocupaci\'on m\'as, y ella puede incitar que usted muestre inseguridad.

La elecci\'on entre un lenguaje riguroso y uno coloquial debe tomar en cuenta la comprensi\'on que desee en la audiencia. Ante un p\'ublico altamente calificado usted puede ahorrarse tiempo al utilizar el lenguaje riguroso, dado que se espera que todos quienes lo escuchan lo dominan. Un lenguaje coloquial suele ser \'util cuando ciertos conceptos y resultados no son conocidos por el p\'ublico; en esos casos puede ser mejor dar breves ideas intuitivas incluyendo nombres f\'aciles de recordar de conceptos y resultados, evitando la rigurosidad que podr\'{\i}a requerir una mayor cantidad de tiempo para su cabal comprensi\'on.

Es dif\'{\i}cil encontrar un lenguaje apropiado para hablar frente a p\'ublico con bajos conocimientos en el tema que se expone.
Imag\'{\i}nese a usted mismo hablando de \'algebra lineal ante personas que no tienen formaci\'on universitaria en temas matem\'aticos. Es probable que el tema del que usted est\'e hablando sea muy adecuado para su audiencia, despu\'es de todo el \'algebra lineal tiene muchas aplicaciones en la vida real. Sin embargo, la audiencia no est\'a calificada para entender t\'erminos como `linealmente independientes' o `que generan un subespacio de dimensi\'on $3$', !`y nosotros no estamos acostumbrados a hablar de estos conceptos con otras palabras! En estas circunstancias debe usted resistirse a su instinto de definir formalmente independencia lineal (piense lo extra\~no que la definici\'on le parecer\'a a quienes le escuchan), y buscar salidas alternativas. En charlas de este estilo se justifica plenamente usar palabras coloquiales para referirse a estos conceptos, y describirlos de manera intuitiva. Por ejemplo, puede aclarar que llamar\'a {\em aplanada} a una terna de vectores cuando se queden en un mismo plano, y sustituir el uso de `linealmente independientes' por `no aplanados'. Note que el t\'ermino ` terna aplanada' es preferible sobre `terna buena' dado que el primero des\-cribe el concepto hasta cierto grado, mientras que el segundo solo hace referencia a una opini\'on subjetiva.


Por otro lado, suponga que entre los asistentes hay alguno que no asocie f\'acilmente los t\'erminos `independencia' o `lineal' con el concepto de `independencia lineal', y ese concepto es definido sin suficiente motivaci\'on ni suficientes ejemplos (es decir, sin dedicarle suficiente tiempo). Una consecuencia frecuente es que esa persona deje de comprender los enunciados en que se menciona independencia lineal, y en esos tramos de la pl\'atica se dedique a verificar las cuentas realizadas, sin entender las ideas detr\'as de ellas. En otras ocasiones quien escucha la pl\'atica da por perdida la esperanza de seguir al expositor hasta el final de la charla, y cambia su atenci\'on hacia sus asuntos personales. Le invito a que compare estos resultados con el objetivo que se plante\'o para su pl\'atica.

Cada \'area de matem\'aticas tiene un lenguaje especializado, que se ha moldeado a lo largo de los a\~nos. Gente que no se dedica a nuestra \'area no necesariamente usa el mismo lenguaje, a\'un cuando use los mismos conceptos. 

Tome usted por ejemplo el t\'ermino `variedad'. Si bien la noci\'on que un top\'ologo, un ge\'ometra diferencial y un ge\'ometra algebraico tengan de estos objetos puede ser compatible, es posible que surjan diferencias cuando se abordan resultados espec\'{\i}ficos. Puede ocurrir que el expositor permita a la variedad tener un conjunto finito de puntos singulares, y que esto est\'e prohibido en el \'area de algunos de la audiencia; o que por default alguien m\'as asuma orientabilidad de las variedades, contrario a la teor\'{\i}a m\'as general que admite variedades no orientables. Esto suele causar confusi\'on (a veces peque\~na y a veces grande) en la audiencia.

Un ejemplo com\'un de lo anterior se da en pl\'aticas de topolog\'{\i}a, en las que se usa la palabra `funci\'on' para referirse exclusivamente a `funciones continuas', pues son las \'unicas relevantes para los resultados de la pl\'atica. Hacer esto frente a personas que trabajan en otras \'areas de matem\'aticas tambi\'en puede originar confusi\'on.
Usted evitar\'a este tipo de confusiones si, o bien elude los t\'erminos que no son est\'andares en su audiencia, o bien anuncia desde el principio de su charla el significado que dar\'a de esos t\'erminos.

Es deseable motivar los temas de los que hablamos resaltando su relevancia hist\'orica o reciente. Sin embargo esto puede requerir cierto tacto. Se debe tener cuidado con los superlativos como `el mayor matem\'atico de la historia' o `el \'area m\'as bonita de matem\'aticas', sobre todo cuando en la audiencia hay personas de \'areas distintas a la nuestra. Estas opiniones son subjetivas y pueden causar suspicacias acerca de qu\'e tan informado est\'e el expositor, o confrontaci\'on de ideas innecesaria para la pl\'atica. Por otro lado, datos bien documentados u opiniones justificadas con citas verificables deben ser siempre bien recibidos.

\section{Consideraciones generales}

Cuando pedimos sugerencias a distintas personas para la preparaci\'on de nuestra siguiente pl\'atica escuchamos respuestas de muchos tipos. Algunas de ellas dan alg\'un consejo pr\'actico, aunque no siempre apropiado para todas las pl\'aticas, en particular para la pr\'oxima. A continuaci\'on enlisto algunos ejemplos.

{\em (a) Toda pl\'atica debe tener una demostraci\'on.}

{\em (b) Inicie la pl\'atica con un chiste.}

{\em (c) Divida la pl\'atica en tres partes; la primera debe ser comprensible para toda su audiencia, la segunda para las personas de su \'area, y la tercera solo para sus colegas m\'as cercanos.}

En algunas ocasiones este tipo de sugerencias ser\'an apropiadas, pero definitivamente no son universales y deber\'a seguirlas \'unicamente cuando crea que convienen a su pl\'atica. Piense que una de las pl\'aticas m\'as famosas y exitosas en la historia de las matem\'aticas (si no es que la m\'as famosa y exitosa) es la pl\'atica de D. Hilbert en la que present\'o los 23 problemas a atacarse durante el siglo XX. La exposici\'on de Hilbert no hubiera sido mejor si hubiera incluido demostraciones, y para el era deseable que la pl\'atica en su totalidad fuera accesible para los matem\'aticos presentes. Compare esa pl\'atica con las sugerencias enlistadas anteriormente.

No siempre contamos con las mismas herramientas para impartir una pl\'atica. En ocasiones se nos solicita que sea en pizarr\'on, mientras que en otras es necesario preparar una presentaci\'on para ser proyectada. Cuando usted tenga la libertad de elegir, busque la opci\'on que considere m\'as adecuada para los objetivos que se ha planteado para su pl\'atica; cada una tiene sus ventajas y desventajas.

Las presentaciones proyectadas aportan informaci\'on sin requerir de tiempo para escribirlas en un pizarr\'on, lo que contribuye a que el contenido pueda presentarse en menor tiempo. Esto hace que las presentaciones con computadora sean id\'oneas para pl\'aticas cortas. Al hacerlas es importante evitar avanzar demasiado r\'apido. Si bien el expositor puede modificar la rapidez en la que se exponen las ideas, no puede controlar la rapidez en la que las procesan los que le escuchan. No debe extra\~nar que en seminarios de una hora haya quienes soliciten espec\'{\i}ficamente que las pl\'aticas sean en pizarr\'on.

Elaborar una buena presentaci\'on para acompa\~nar una pl\'atica no es tarea f\'acil. La presentaci\'on no debe ser un acorde\'on de donde nosotros podamos leer cada palabra que hemos de decir; no debe ser hecha para ayudar al expositor sino para facilitar la comprensi\'on de la audiencia. Es muy importante entender que las palabras que vayan a aparecer proyectadas no necesariamente van a ser le\'{\i}das, sobre todo si el expositor habla durante el tiempo en que se proyecta la l\'amina que contiene esas pa\-la\-bras. El texto que proyectemos debe ser complemento de lo que digamos. Un texto largo invita a no leerlo y a esperar a que el expositor explique; si el texto no va a ser le\'{\i}do habr\'a que preguntarse la raz\'on por la que fue escrito. Por ello se sugiere no basar su presentaci\'on en l\'aminas que muestren p\'aginas de art\'{\i}culos, o en l\'aminas que contengan muchas letras. En cambio, se recomienda incluir diagramas o dibujos que ilustren de lo que se habla. Tambi\'en se pueden incluir a manera de esquema solo las ideas centrales de lo que se dice durante el tiempo en que se muestra la l\'amina.


\section{Tome esto en cuenta durante su pl\'atica.}

\begin{wrapfigure}{r}{0.45\textwidth}
\centering
\includegraphics[width=0.45\textwidth]{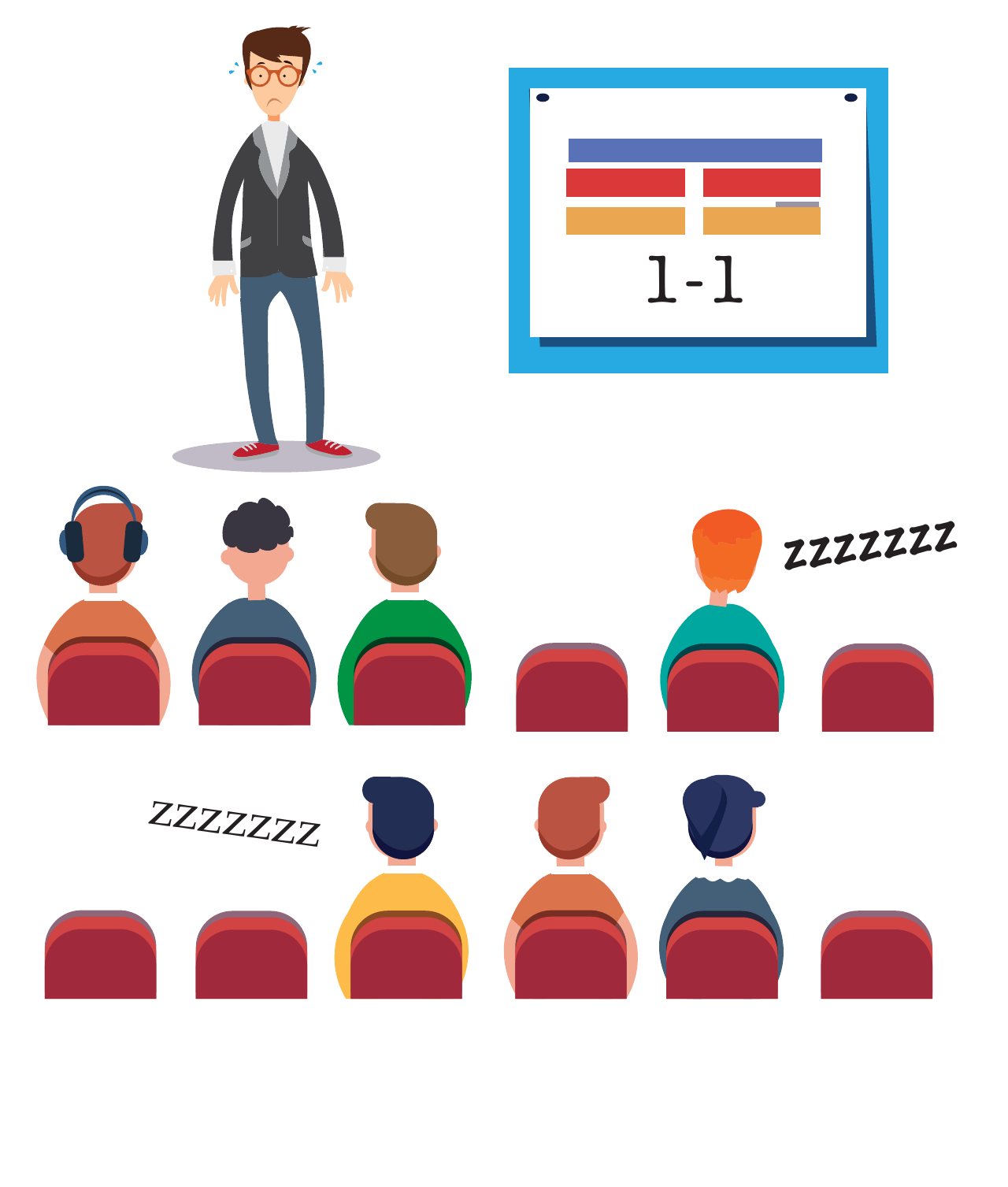}
 \end{wrapfigure}

Es normal que estemos nerviosos durante la pl\'atica, sobre todo en nuestras primeras exposiciones. Hay factores que contribuyen a ello, y tambi\'en maneras de atacarlos.

En primer lugar, el tono de voz y la postura corporal pueden influir fuertemente a crear una atm\'osfera adecuada o adversa. Alguien que habla al pizarr\'on, que no ve a la audiencia y que usa un tono de voz muy bajo invita a dejar de prestar atenci\'on. 
Darse cuenta de una baja en la atenci\'on que se le brinda contribuir\'a a poner m\'as nervioso al expositor.
Debido a lo anterior es conveniente buscar contacto visual constante con la audiencia y modular el tono de voz de modo que se resalten las cosas m\'as importantes (se puede subir un poco el tono en ellas), y nunca se hable demasiado bajo como para que no se escuche bien lo que decimos.

Cuando estamos nerviosos somos susceptibles a crear impresiones err\'oneas en nosotros mismos, en particular al tratar de interpretar las caras de quienes nos escuchan. Una de estas impresiones es asumir que en la audiencia hay personas que se est\'an aburriendo, y podemos cometer el error de incrementar la velocidad en la que abordamos los temas. Otra m\'as es cuando las caras de algunos asistentes nos sugieren que no est\'an entendiendo, y ponemos m\'as atenci\'on a explicar con m\'as detalle de lo planeado. Debemos evitar tomar estas impresiones como hechos consumados; de ser posible, debemos ignorar impresiones que no tengan bases s\'olidas. Piense que es posible que, a pesar de las apariencias, a la audiencia le estaba costando trabajo seguir la pl\'atica, y acelerar el ritmo solo empeora la situaci\'on. Por otro lado, abundar en las explicaciones m\'as de lo planeado suele romper el flujo l\'ogico de las ideas. De no tener evidencia convincente de que algo est\'a saliendo mal, !`contin\'ue con lo planeado!

A todos nos sucede de vez en cuando que hay un error en nuestras l\'aminas y nos damos cuenta a la hora de la presentaci\'on. Hay varias maneras de abordarlo. Aqu\'{\i} le brindo unos ejemplos.\newline
{\em Seria y formal:} aclarar el error y pedir disculpas de manera breve; poste\-rior\-mente continuar como si nada hubiera ocurrido.\newline
{\em Bromista:} enfatizar en la presunta veracidad del texto con el error (no m\'as de unos segundos), para luego aclarar lo que deb\'{\i}a haberse escrito.\newline
{\em Indiferente:} evitar leer el texto con error y hacer como si no estuviera ah\'{\i}.\newline
Cada quien debe elegir c\'omo proceder de acuerdo a su car\'acter, y tener en cuenta que la salida que elijamos debe cumplir con que
\begin{itemize}
 \item no se da informaci\'on falsa,
 \item no se utiliza mucho tiempo en remediar la situaci\'on,
 \item usted no se distrae por el error y retoma el flujo de la pl\'atica cuanto antes.
\end{itemize}
Hay algunas ocasiones en que se deben evitar algunas de estas estrategias. Piense que si un error es evidente, eludir su menci\'on puede ser una estrategia equivocada, pues es posible que se le interrumpa para verificar que se trata de un error. Por otro lado, un error ortogr\'afico dif\'{\i}cilmente causa confusi\'on, y no merece ocupar tiempo alguno en evidenciarlo. 

En secciones anteriores suger\'{\i} evitar ser formal en partes de la pl\'atica, ya sea por cuestiones de tiempo, de fluidez de las ideas o de probabilidad de comprensi\'on por parte de la audiencia. Esto no debe interprestarse como una invitaci\'on a abandonar la formalidad por completo, sobre todo si el omitir detalles deja una impresi\'on falsa del tema en cuesti\'on en quienes le escuchan.

Por sobre todas las cosas le sugiero evitar decir mentiras a su audiencia. Aun si ignora la respuesta a alguna pregunta que le formulen, !`Es mejor decir `no s\'e' a decir mentiras! Tenga en cuenta que algunos comentarios y preguntas en las pl\'aticas tienen como objetivo sugerir al expositor direcciones hacia d\'onde continuar, y no se espera que sepa a detalle todo acerca de esas direcciones. La intenci\'on de quienes preguntan en pl\'aticas casi nunca es evidenciar que hay algo que el expositor no sabe.

\section{Despu\'es de la pl\'atica a\'un hay cosas que usted puede hacer}

Hay aspectos sencillos de evaluar de la pl\'atica que reci\'en impartimos. Uno de ellos es el tiempo; ?`lo planemos bien? ?`nos falt\'o? ?`nos sobr\'o?

Detectar otros aspectos requiere de voluntad autocr\'{\i}tica. Cuesta trabajo darnos cuenta de si hablamos demasiado r\'apido o en un volumen demasiado bajo. Para ese tipo de cosas es bueno contar con alguien de confianza entre los asistentes que nos pueda hacer ver esos detalles.

La misma audiencia aporta evidencia de qu\'e tan bien o qu\'e tan mal estuvo la pl\'atica. Interrupciones pertinentes suelen indicar que esas personas est\'an preocupadas por seguir lo que decimos. Solicitudes de ejemplos, o de volver a enunciar un resultado anterior indican que el tema y la manera de exponer captaron la atenci\'on de esas personas. Una pl\'atica sin interrupciones y cuya \'unica pregunta al final busca relacionar lo expuesto con otros temas de matem\'aticas es un indicio 
de que no fueron muchos los que acompa\~naron al expositor hasta el final. Repetidas solicitudes de definiciones de conceptos b\'asicos son clara muestra de que la pl\'atica fue de nivel muy elevado para la audiencia a la que fue presentada.

Cuando damos una pl\'atica desastrosa tendemos a desanimarnos y suele ocurr\'{\i}rsenos dar el menor n\'umero de pl\'aticas posibles en el futuro. Sin embargo, lo que debemos hacer !`es precisamente lo contrario! Primero hay que identificar aspectos negativos de nuestra pl\'atica que queremos que no se repitan, y de ser posible tambi\'en aspectos positivos que queremos que sigan presentes en nuestras exposiciones. Despu\'es hay que buscar foros adecuados para dar pl\'aticas, como seminarios formales o informales. En cada pl\'atica hay que buscar corregir alg\'un aspecto negativo.

Es posible que para dar una pl\'atica que finalmente le deje satisfecho usted deba impartir otras cinco ponencias intermedias; si persevera, tenga por seguro que llegar\'a ese d\'{\i}a. En cambio, si decide dar el menor n\'umero de pl\'aticas posibles puede ser que pase una vida completa en la que dar pl\'aticas represente angustias y desilusiones. !\`Animo! !`Todos podemos dar una buena pl\'atica!




\section{Acerca de este texto}

El presente texto es una versi\'on preliminar de un trabajo enviado a Miscel\'anea Matem\'atica. Fue escrito en M\'exico en 2017 y es producto de la preocupaci\'on de varios colegas por la baja calidad de muchas pl\'aticas de matem\'aticas que presenciamos en seminarios, coloquios y congresos. Los s\'{\i}ntomas de este hecho son de diversas naturalezas y algunos se describen a continuaci\'on.

Hoy en d\'{\i}a la asistencia a nivel nacional a coloquios y seminarios de becarios de matem\'aticas no es lo nutrida que deber\'{\i}a. Cuando se cuestiona de esto a estudiantes y acad\'emicos se deja entrever que para muchos de ellos es una p\'erdida de tiempo asistir a estos eventos. Una parte fundamental de esta opini\'on es que el material presentado les resulta incomprensible.

Vemos con preocupaci\'on que asistentes a pl\'aticas, cada vez con mayor frecuencia, pasan una buena parte de la exposici\'on frente al monitor de su computadora port\'atil, una tablet o un celular. No es raro ver que en esta situaci\'on alguna pantalla muestra una p\'agina de alguna red social o un videojuego.

Es com\'un que haya pocas preguntas al final de las charlas, y que sean todav\'{\i}a menos las que son relevantes al tema expuesto. Algunas son por cortes\'{\i}a y son formuladas por alg\'un experto en el \'area o por el moderador con la intenci\'on de que el expositor no se desanime.

Adem\'as de los s\'{\i}ntomas antes mencionados, notamos falta de conciencia de la comunidad matem\'atica de la gravedad de esta situaci\'on, y por lo mismo existen pocas acciones orientadas a remediarla.

Por si fuera poco, hay una presi\'on constante a dar pl\'aticas para fines de apoyos econ\'omicos, obtenci\'on de grado, contrataciones y promociones; pero no incluye una solicitud de calidad en estas pl\'aticas. En comparaci\'on, la presi\'on para que la gente escuche pl\'aticas es casi nula. De acuerdo a lo anterior, el futuro que nos espera es que las pl\'aticas no sean entendidas por nadie de los que las escuchan, o peor aun, que no haya quien las escuche. Las agencias que dan apoyos financieros y los comit\'es que deciden obtenciones de grado, contrataciones y promociones no notar\'{\i}an si alguno de estos fuera el caso.

Este texto busca difundir la existencia del problema y solicitar la parti\-ci\-paci\'on de todos sus lectores para colaborar en posibles soluciones. Ser\'a un gran paso que cada lector decida no ser parte del problema y encargarse de dar pl\'aticas de buena calidad.

El presente texto no trata de decir que hay una receta para dar una pl\'atica perfecta. El autor de estas l\'{\i}neas est\'a convencido de que cada expositor le debe poner un toque de su personalidad a cada una de sus pl\'aticas, sin importar si se trata de una persona preferentemente seria o bromista. Lo central en las mejores pl\'aticas que he escuchado no es el tono del expositor, sino una serie de factores que hacen que me halle cautivo a lo largo de toda la exposici\'on.

No espero que el total de la comunidad matem\'atica concuerde conmigo en cada uno de los temas que abordo en este texto.
Por el contrario, atentamente invito a todo matem\'atico a asistir a pl\'aticas y ponerles la mayor atenci\'on posible. Al hacerlo sugiero realizar el ejercicio de determinar qu\'e elementos le agradaron y cu\'ales no de modo que busque adoptar los primeros en la medida de lo posible y evitar caer en los segundos, independientemente de si ello contradice lo escrito aqu\'{\i}.

Tambi\'en le recomiendo al lector buscar otros textos cuyo prop\'osito sea orientar en la preparaci\'on de pl\'aticas con el objetivo de tener m\'as puntos de comparaci\'on. Algunos de ellos, con distintos objetivos particulares, se enlistan a continuaci\'on. Comenzamos con \cite{Halmos}, un texto escrito hace m\'as de 40 a\~nos a solicitud de una sociedad matem\'atica enfocado a presentaciones ante matem\'aticos no especializados en el \'area.
El texto \cite{Matilde} es la base de un curso orientado a expresarse bien al dar charlas, especialmente ante audiencias dif\'{\i}ciles.
En \cite{Agelos} se muestra una preocupaci\'on similar a la presentada aqu\'{\i}, aunque enfatiza otros aspectos.
Finalmente, los textos \cite{Continuos1} y \cite{Continuos2} muestran lineamientos y requisitos para impartir pl\'aticas de 20 minutos en un evento anual de un \'area espec\'{\i}fica de topolog\'{\i}a


Finalmente, pido disculpas por anticipado si alg\'un lector se sintiera ofendido por los ejemplos que menciono de elementos negativos de una pl\'atica. Ninguno ha sido incluido con la intenci\'on de atacar a alg\'un \'area de matem\'aticas ni a ning\'un expositor en particular.

\section{Agradecimientos}
El autor agradece a Coppelia Cerda por la elaboraci\'on de las ilustraciones, a Amanda Montejano, Patricia Pellicer, Ferr\'an Valdez y Jos\'e Antonio Montero y al r\'eferi anónimo cuyas sugerencias contribuyeron a hacer de este un mejor texto.

\bibliographystyle{miscelanea}
\bibliography{refer}

\end{document}